\documentclass{article}
\usepackage{PRIMEarxiv}
\usepackage{amsmath}
\usepackage{subcaption}
\usepackage[linesnumbered,ruled,vlined]{algorithm2e}
\usepackage{float}
\usepackage[utf8]{inputenc}
\usepackage[T1]{fontenc}
\usepackage{hyperref}
\usepackage{url}
\usepackage{booktabs}
\usepackage{amsfonts}
\usepackage{nicefrac}
\usepackage{microtype}
\usepackage{lipsum}
\usepackage{fancyhdr}
\usepackage{graphicx}
\usepackage{tikz}
\usepackage{pgfplots}
\pgfplotsset{compat=1.18}
\usetikzlibrary{arrows.meta, positioning, shapes.geometric, calc, fit}
\graphicspath{{media/}{figuras/}}

\title{The Birthday Problem and the Birthday Person Problem: Intuition, Probability and Cryptographic Applications.}

\author{
  Alisson Garcia da Silva\\
  Federal University of Rio Grande do Norte \\
  Caicó- RN\\
  \texttt{alisson.silva.127@ufrn.edu.br} \\
   \And
   Iarythssa Duarte de Araújo\\
   Department of Exact and Applied Sciences  \\
   Federal University of Rio Grande do Norte \\
   Caicó- RN\\
  \texttt{m.sc.iarythssa@gmail.com} \\
   \And
  Francisco Márcio Barboza \\
  Department of Computing and Technology \\
  Federal University of Rio Grande do Norte \\
  Caicó- RN\\
  \texttt{marcio.barboza@ufrn.br} \\
}

\begin{document}
\maketitle

\begin{abstract}
The present work discusses the birthday problem and the birthday person problem, two classical problems in probability theory widely used to illustrate situations in which human intuition diverges from mathematical results. Taking the work of \cite{magalhaes2006} as its main reference, the study analyzes the probability of birthday coincidences within groups of people, considering a year with \(365\) days and assuming that all birth dates are equally likely. The birthday problem seeks to determine the probability that at least two people in a group share the same birthday, showing that, even in relatively small groups, the probability of a coincidence increases rapidly due to the large number of possible pairwise comparisons among individuals. In contrast, the birthday person problem investigates the probability that a previously selected person shares the same birthday as another member of the group. This probability increases more slowly, since the number of relevant comparisons is restricted to the pairs formed by the selected person and the remaining members of the group. Despite the interpretative differences, both problems share a similar mathematical foundation based on combinatorial analysis and probability theory. In addition to the theoretical development, numerical results and graphical representations are presented to illustrate the behavior of the probabilities involved and to visualize how the likelihood of coincidences increases as the number of people grows. Finally, applications related to computing and cryptography are discussed, particularly the study of collisions in hash functions and the so-called birthday attack. The discussion highlights the importance of probabilistic reasoning in understanding apparently unlikely phenomena, while also emphasizing how the growth in the number of possible combinations directly influences the results obtained.

\end{abstract}


\section{Introduction}

The study of coincidences is a recurring topic in probability theory and combinatorial analysis, particularly in problems involving independent and identically distributed random events. Among the best-known examples of this type is the birthday problem, whose formulation investigates the probability that, in a group of $n$ people, at least two of them share the same birthday. Despite its simple formulation, the problem has become one of the classical examples in probability theory because of the apparently unexpected nature of the results obtained.

Historically, the birthday problem has attracted considerable attention because it highlights the discrepancy between intuitive reasoning and probabilistic behavior. As discussed by \cite{diaconis1989}, the birthday problem serves as an important starting point for the probabilistic study of coincidences, particularly in situations involving surprise and the subjective interpretation of results.

In general, it is common to imagine that birthday coincidences are rare in small groups of people. However, probabilistic results show a behavior different from that initially suggested by intuitive perception. In a group with only \(23\) people, for example, the probability that at least two of them share the same birthday already exceeds $50\%$ \cite{ross2019a}. This result is often associated with the so-called ``birthday paradox'', not because it represents a logical contradiction, but because it highlights the contrast between probabilistic calculation and intuitive interpretation.

The seemingly paradoxical character of the problem is directly related to the growth of the number of distinct pairs of individuals that can be formed in a group. While intuition often considers comparisons that are individual or related to only one fixed person, the classical probabilistic model simultaneously involves all possible comparisons between pairs of individuals. In this way, the total number of comparisons grows rapidly as the number of people considered increases, being described by the expression

\[
\binom{n}{2}=\frac{n(n-1)}{2},
\]

which explains the accelerated growth of coincidence probabilities even in relatively small groups.

The mathematical development of the birthday problem is based on classical concepts of probability theory and combinatorial analysis, especially those presented by \cite{magalhaes2006}. Among the main concepts used are sample space, events, event independence, complementary events and the multiplicative principle, which provide the necessary theoretical basis for the formulation of the problem.

In addition to the classical formulation, the literature presents other variations of the problem, among them the so-called Birthday Person problem. In this version, the probability of a previously fixed individual sharing his or her birthday with at least one of the other members of the group is investigated. As discussed by \cite{amato2020}, different versions of the problem can be interpreted using combinatorial and probabilistic techniques, enabling associations with more general mathematical structures such as graphs and chromatic polynomials.

Currently, the birthday problem goes beyond the purely theoretical field of probability, having applications in areas related to computing, cryptography, information security and hash functions. In particular, probabilistic concepts associated with the problem are used in the analysis of collisions in hash functions and in the study of the so-called birthday attack, phenomena that play an important role in assessing the security of cryptographic systems \cite{stallings2015}.

In this way, the birthday problem stands out not only as a classic example of probability theory, but also as an important mathematical model capable of relating probability, combinatorial analysis and computational applications. In this context, the present work aims to develop mathematically the birthday problem and the birthday person problem, analyzing the growth of coincidence probabilities, their combinatorial interpretations and some applications associated with the problem.

\section{Theoretical Background}

\subsection{Basic Concepts of Probability}

The development of the birthday problem is based on classical concepts of probability theory and combinatorial analysis, especially those related to the study of random events, event independence and complementary probabilities. In this context, the concepts of sample space, events, disjoint events and Kolmogorov's axioms stand out, which provide the necessary mathematical basis for the probabilistic formulation of the problem.

The sample space is defined as the set of all possible outcomes associated with a random experiment, and is usually represented by \(\Omega\). According to \cite{magalhaes2006}, the subsets of the sample space are called events, representing specific situations associated with the experiment performed. Thus, in the context of the birthday problem, an event can represent, for example, the occurrence of a coincidence between birthday dates in a group of people.

Among the fundamental concepts used in probability theory, the complementary event stands out. Given an event \(A\), its complement, represented by \(A^c\), corresponds to the set of all outcomes in the sample space that do not belong to \(A\). In symbols,

\[
A^c=\Omega\setminus A.
\]

In the birthday problem, the calculation of the desired probability is often carried out by means of complementary probability, that is, by first calculating the probability that no birthday coincidences occur.

Another important concept is that of disjoint events. Two events \(A\) and \(B\) are called disjoint when they cannot occur simultaneously, that is,

\[
A \cap B = \emptyset.
\]

The classical definition of probability considers that, in experiments with equiprobable outcomes, the probability of occurrence of an event \(A\) is given by the ratio of the number of favorable cases to the total number of possible cases; that is,

\[
P(A)=\frac{n(A)}{n(\Omega)},
\]

where \(n(A)\) represents the number of outcomes favorable to the event and \(n(\Omega)\) the total number of elements in the sample space.

\subsection{Definition 1.1 ($\sigma$-algebra)}

A class of subsets of $\Omega$, represented by $\mathcal{F}$, is called a $\sigma$-algebra if it satisfies the following properties:

\begin{itemize}
    \item[(A1)] $\Omega \in \mathcal{F}$;
    \item[(A2)] If $A \in \mathcal{F}$, then $A^{c} \in \mathcal{F}$;
    \item[(A3)] If $A_i \in \mathcal{F}$ for $i \geq 1$, then
    \[
    \bigcup_{i=1}^{\infty} A_i \in \mathcal{F}.
    \]
\end{itemize}

In addition to the classical definition, the modern theory of probability is based on the axiomatic formulation proposed by Kolmogorov \cite{magalhaes2006}. According to this approach, a probability must satisfy three fundamental axioms.

\subsection{Definition 1.2 (Probability)}

A function \(P\), defined on the \(\sigma\)-algebra \(\mathcal{F}\) of subsets of \(\Omega\) and taking values in \([0,1]\), is called a probability if it satisfies Kolmogorov's axioms:

\begin{itemize}
    \item[(Ax1)] \(P(\Omega)=1\);
    \item[(Ax2)] For every \(A\in\mathcal{F}\), \(P(A)\geq 0\);
    \item[(Ax3)] For every sequence of mutually exclusive events \(A_1,A_2,\ldots\in\mathcal{F}\), we have
    \[
    P\left(\bigcup_{i=1}^{\infty}A_i\right)
    =
    \sum_{i=1}^{\infty}P(A_i).
    \]
\end{itemize}

\subsection{Independence of Events}

In probability theory, a fundamental concept is that of independence between events. According to \cite{magalhaes2006}, two events are independent when the occurrence of one of them does not interfere with the probability of the occurrence of the other.

\subsection{Definition 1.4 (Independence of two events)}

Let \(A\) and \(B\) be events of a probability space $(\Omega,\mathcal{F}, P)$. We say that \(A\) and \(B\) are independent if the information of the occurrence of \(B\) does not change the probability assigned to the event \(A\), that is,

\[
P(A\cap B)=P(A)P(B).
\]

When \(P(B)>0\), this condition is equivalent to

\[
P(A\mid B)=P(A).
\]

Similarly, when \(P(A)>0\),

\[
P(B\mid A)=P(B).
\]

This equality shows that the probability of the simultaneous occurrence of the events can be obtained by the product of their individual probabilities. From a probabilistic point of view, independence indicates that knowledge of the occurrence of one event does not modify the probabilistic behavior of the other. Thus, the occurrence of one event does not provide additional information about the occurrence of the second event.

In the birthday problem, it is assumed that people's birthdays are independent of each other. This means that one person's birthday does not interfere with the probabilistic distribution associated with the birthdays of the other people in the group.

This hypothesis is essential for the development of the probabilistic expressions used in the problem. By calculating the probability that several people have distinct birthdays, it becomes possible to successively multiply the probabilities associated with each stage of the experiment.

Therefore, event independence constitutes one of the main probabilistic foundations used in the mathematical formulation of both the birthday problem and the birthday person problem.

\section{The Birthday Problem}

\subsection{Formulation of the problem}

Determine the probability that, in a group with \(n\) people, at least two of them have the same birthday.

\subsection{Solution to the problem}

Initially, for the probabilistic modeling of the problem, we adopt the following hypotheses: a year with \(365\) days, in which all dates are equally probable; birthdays are independent; and February \(29\) is disregarded.

Let the event be

\[
A=\{\text{at least two people have the same birthday}\}.
\]

We want to calculate the probability of occurrence of the event \(A\), that is,

\[
P(A).
\]

To simplify the calculations, we initially consider the complementary event of \(A\), defined by

\[
A^c = \{\text{all people have distinct birthdays}\}.
\]

Note that, since \(A\) and \(A^c\) are complementary events in \(\Omega\),

\[
A\cup A^c=\Omega.
\]

Necessarily, one of the two events occurs in our sample space: either there is a birthday coincidence, or all birthday dates are distinct. Furthermore,

\[
A\cap A^c=\emptyset,
\]

since an event and its complement cannot occur simultaneously. Since \(A\) and \(A^c\) are disjoint, by the additivity axiom,

\[
P(A\cup A^c)=P(A)+P(A^c).
\]

But

\[
A\cup A^c=\Omega.
\]

We obtain

\[
P(\Omega)=P(A)+P(A^c).
\]

By Kolmogorov's first axiom,

\[
P(\Omega)=1.
\]

Substituting \(P(\Omega)=1\), we obtain

\[
1=P(A)+P(A^c).
\]

Now, isolating \(P(A)\), one subtracts \(P(A^c)\) from both sides of the equality,

\[
1-P(A^c)=P(A)+P(A^c)-P(A^c).
\]

Since

\[
P(A^c)-P(A^c)=0,
\]

it follows that

\[
1-P(A^c)=P(A).
\]

Therefore,

\[
P(A)=1-P(A^c).
\]

Thus, it is sufficient to calculate \(P(A^c)\), that is, to determine the probability that all people have distinct birthdays.

Now we calculate this probability. For this, we define the following events:

\begin{itemize}
    \item \(E_1\): the first person has any birthday;
    \item \(E_2\): the second person has a birthday different from the first;
    \item \(E_3\): the third person has a birthday different from the previous two;
    \item \(\cdots\)
    \item \(E_n\): the \(n\)-th person has a birthday different from all previous ones.
\end{itemize}

Thus, under these conditions,

\[
A^c=E_1\cap E_2\cap E_3\cap\cdots\cap E_n.
\]

Note that the events \(E_1,E_2,E_3,\ldots,E_n\) are not independent of each other, since the fact that some people already have distinct birthdays changes the restrictions imposed on the next people in the group. Thus, by the product rule for conditional probabilities,

\[
P(A^c)
=
P(E_1)\,
P(E_2\mid E_1)\,
P(E_3\mid E_1\cap E_2)
\cdots
P(E_n\mid E_1\cap E_2\cap\cdots\cap E_{n-1}).
\]

The first person can be born on any of the \(365\) days of the year. Thus,

\[
P(E_1)=\frac{365}{365}.
\]

The second person must be born on a different day from the first. Therefore,

\[
P(E_2\mid E_1)=\frac{364}{365}.
\]

The third person must have a birthday different from those of the previous two people. Thus,

\[
P(E_3\mid E_1\cap E_2)=\frac{363}{365}.
\]

Continuing the process, we obtain

\[
P(A^c)
=
\frac{365}{365}\cdot
\frac{364}{365}\cdot
\frac{363}{365}\cdots
\frac{365-(n-1)}{365}.
\]

Or, equivalently,

\begin{equation}
P(A^c)
=
\frac{365\cdot364\cdot363\cdots(365-n+1)}{365^n}.
\label{eq:complementar}
\end{equation}

We know that

\[
P(A)=1-P(A^c).
\]

Substituting Equation \ref{eq:complementar}, we obtain

\[
P(A)
=
1-
\frac{365\cdot364\cdot363\cdots(365-n+1)}{365^n}.
\]

Therefore, the probability that at least two people have the same birthday is

\begin{equation}
P(A)
=
1-
\frac{365\cdot364\cdot363\cdots(365-n+1)}{365^n}.
\label{eq:probabilidade}
\end{equation}

\subsection{Example 1}
\label{subsec:example1}

Consider a group with \(n=20\) people. We wish to determine the probability that at least two of them have the same birthday. Thus, by Equation \ref{eq:probabilidade}, the probability is given by

\begin{equation}
P(A)
=
1-
\frac{365\cdot364\cdot363\cdots346}{365^{20}}.
\label{eq:exemplo20}
\end{equation}

Initially, note that

\[
P(A^c)
=
\frac{365}{365}\cdot
\frac{364}{365}\cdot
\frac{363}{365}\cdots
\frac{346}{365}.
\]

Performing the successive products, we obtain

\[
P(A^c)\approx0{,}58856.
\]

Now, substituting into Equation \ref{eq:exemplo20}, we obtain

\[
P(A)=1-0{,}58856.
\]

Thus,

\[
P(A)\approx0{,}41144.
\]

Therefore,

\[
P(A)\approx41{,}144\%.
\]

The result indicates that the probability that at least two people have the same birthday in a group of \(20\) people is approximately \(41{,}144\%\).

The result obtained in Subsection \ref{subsec:example1} may seem surprising at first glance, since it is common to imagine that birthday coincidences are rare events in relatively small groups. However, probabilistic analysis shows that the occurrence of these coincidences becomes significantly more likely as the number of possible comparisons between individuals in the group increases.

The seemingly paradoxical behavior of the birthday problem is directly related to the growth of the number of possible comparisons carried out between individuals in a group. In an initial intuitive analysis, it is common to consider only comparisons involving a single person. However, in the classical probabilistic model of the problem, all people in the group are compared with each other simultaneously.

From the combinatorial point of view, the total number of possible comparisons corresponds to the number of distinct pairs that can be formed in a group with \(n\) people. This value is given by

\[
\binom{n}{2}=\frac{n(n-1)}{2},
\]

an expression obtained from the combination of \(n\) elements taken two by two.

The analysis notes that the growth of this expression occurs quadratically with respect to the number of individuals in the group. Thus, small variations in \(n\) produce significant increases in the number of comparisons performed. This behavior explains why coincidence probabilities grow rapidly, even in relatively small groups.

For example, in a group with \(23\) people, the total number of possible pairs is given by

\[
\frac{23(23-1)}{2}=253.
\]

Thus, although the group has only \(23\) individuals, there are \(253\) possible comparisons between pairs of people. This accelerated growth in the number of comparisons constitutes one of the main factors responsible for the behavior often associated with the so-called birthday paradox.

Therefore, the combinatorial interpretation of the problem makes it possible to understand that the growth of the probabilities of coincidence does not occur linearly, but is directly associated with the increase in the number of possible pairs formed in a group.

\subsection{Asymptotic approximation}
\label{subsec:poisson}

For large values of \(n\), the exact formula \eqref{eq:probabilidade} becomes numerically cumbersome. A widely used approximation follows from the fact that, for $x$ small, $\ln(1-x)\approx -x$. Indeed,

\[
P(A^c)=\prod_{k=0}^{n-1}\left(1-\frac{k}{365}\right)
\approx\exp\left(-\frac{1}{365}\sum_{k=0}^{n-1}k\right)
=\exp\left(-\frac{n(n-1)}{2\cdot 365}\right).
\]

Hence,

\begin{equation}
P(A)\approx 1-\exp\left(-\frac{n(n-1)}{730}\right).
\label{eq:poisson}
\end{equation}

Approximation \eqref{eq:poisson} is known as the \emph{Poisson approximation} of the birthday problem and is particularly useful in cryptographic contexts. Solving $P(A)\approx 1/2$ in \eqref{eq:poisson} yields $n\approx 1{,}1774\sqrt{365}\approx 22{,}49$, which corroborates the classical result $n=23$.

\section{Numerical Results}

\subsection{Numerical values of the probabilities}

To analyze the growth of coincidence probabilities in the birthday problem, some numerical results obtained from formula \eqref{eq:probabilidade} are presented. The values in Table \ref{table:birthday_probabilities} make it possible to observe the rapid growth of probabilities as the number of people in the group increases.

\begin{table}[H]
\centering
\begin{tabular}{|c|c|}
\hline
\textbf{Number of people} & \textbf{Probability of coincidence} \\
\hline
\(10\) & \(11{,}69\%\) \\
\(20\) & \(41{,}14\%\) \\
\(23\) & \(50{,}73\%\) \\
\(30\) & \(70{,}63\%\) \\
\(50\) & \(97{,}04\%\) \\
\(60\) & \(99{,}41\%\) \\
\hline
\end{tabular}
\caption{Probability of birthday coincidence according to the number of people in the group.}
\label{table:birthday_probabilities}
\end{table}

\subsection{Analysis of the results}

It can be observed that the growth of the probabilities occurs quite rapidly, as shown in Table \ref{table:birthday_probabilities}. Although, intuitively, it is common to imagine that birthday coincidences are rare, the results show that relatively small groups already produce high probabilities of repeated dates.

In particular, the case of \(23\) people stands out, in which the probability of coincidence exceeds \(50\%\). This behavior is directly related to the rapid growth in the number of possible comparisons between pairs of individuals, a fundamental aspect for understanding the so-called birthday paradox.

Furthermore, it can be observed that the growth of the probabilities does not occur linearly. As the number of people increases, the probabilities rapidly approach \(100\%\), demonstrating the strong influence of the combinatorial growth associated with the problem.

\subsection{Graphical representations}

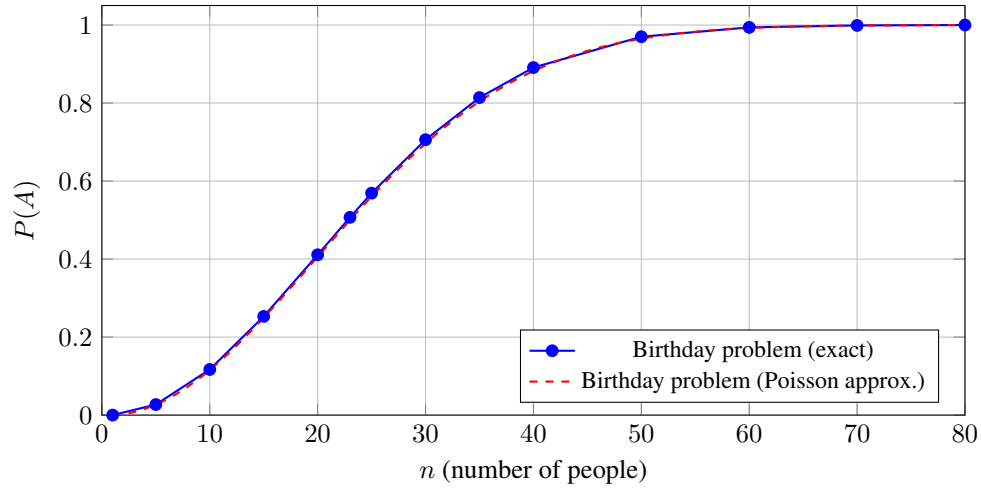
\begin{figure}[H]
\centering
\begin{tikzpicture}
\begin{axis}[
    width=13cm, height=7cm,
    xlabel={$n$ (number of people)},
    ylabel={$P(A)$},
    xmin=0, xmax=80,
    ymin=0, ymax=1.05,
    xtick={0,10,...,80},
    ytick={0,0.2,...,1},
    grid=major,
    legend pos=south east,
    legend style={font=\small},
]
\addplot[blue, thick, mark=*] coordinates {
    (1, 0.000) (5, 0.027) (10, 0.117) (15, 0.253) (20, 0.411)
    (23, 0.507) (25, 0.569) (30, 0.706) (35, 0.814) (40, 0.891)
    (50, 0.970) (60, 0.994) (70, 0.999) (80, 1.000)
};
\addlegendentry{Birthday problem (exact)}
\addplot[red, dashed, thick, domain=1:80, samples=200] {1-exp(-x*(x-1)/730)};
\addlegendentry{Birthday problem (Poisson approx.)}
\end{axis}
\end{tikzpicture}
\caption{Probability of at least one birthday coincidence as a function of the group size \(n\), comparing the exact formula \eqref{eq:probabilidade} with the Poisson approximation \eqref{eq:poisson}.}
\label{fig:curva_aniversarios}
\end{figure}

\begin{figure}[H]
    \centering
    \includegraphics[width=15cm]{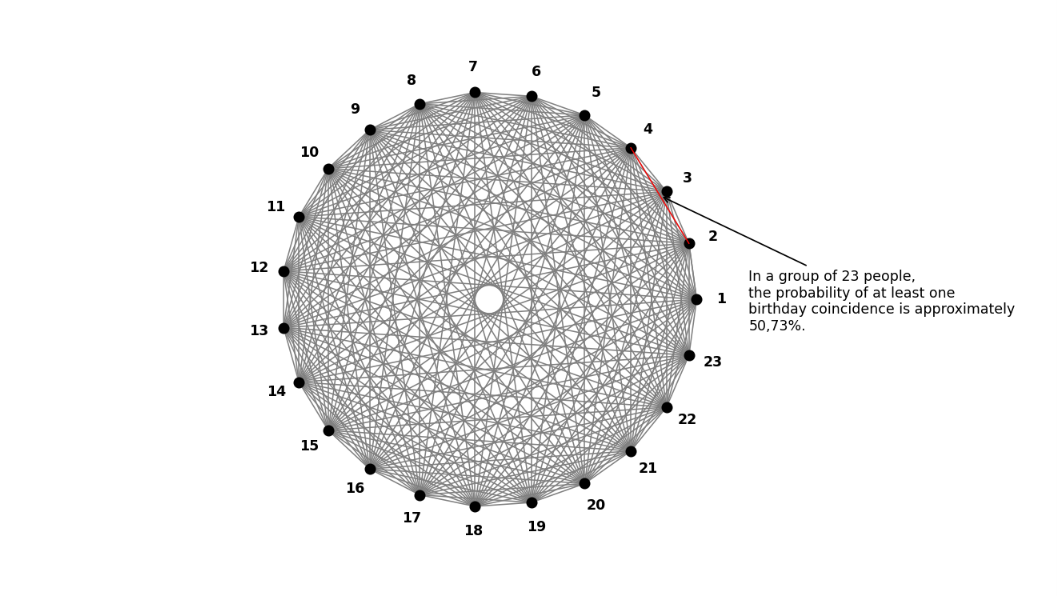}
    \caption{Schematic representation of the birthday problem for a group of 23 people. Each vertex represents an individual and each edge a pairwise comparison; the red edge highlights one of the $253$ possible coincidences.}
    \label{fig:aniversarios23}
\end{figure}

\begin{figure}[H]
    \centering
    \includegraphics[width=15cm]{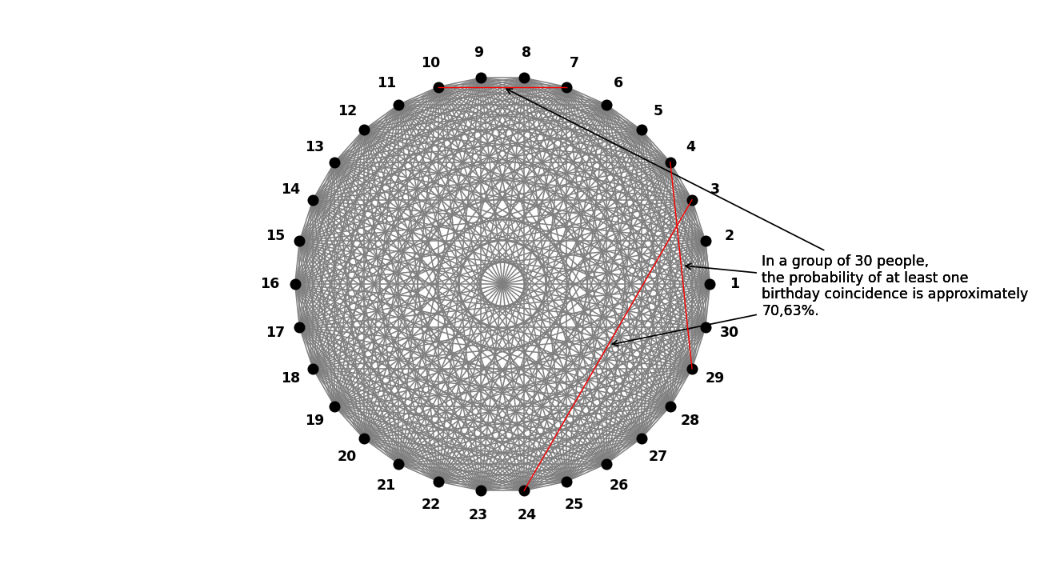}
    \caption{Schematic representation of the birthday problem for a group of 30 people. The number of pairwise comparisons grows quadratically ($\binom{30}{2}=435$), which explains the high coincidence probability.}
    \label{fig:aniversarios30}
\end{figure}

To complement the theoretical analysis and numerical results presented, graphical representations are provided. Figure \ref{fig:curva_aniversarios} shows the exact probability curve together with its Poisson approximation \eqref{eq:poisson}, allowing one to visualize the excellent adherence of the approximation even for moderate values of \(n\). Figures \ref{fig:aniversarios23} and \ref{fig:aniversarios30} present schematic representations of the birthday problem for groups of \(23\) and \(30\) people, respectively, in which each vertex corresponds to an individual and each edge represents a pairwise comparison. The number of edges grows quadratically with \(n\), illustrating the combinatorial explosion responsible for the high coincidence probabilities.

In this way, the visual resources contribute to a better understanding of the probabilistic behavior analyzed, making the relationship between the calculations developed and the phenomena observed in practice clearer.

\section{The Birthday Person Problem}

In addition to the classical formulation of the birthday problem, there is also the so-called Birthday Person problem, in which comparisons are carried out with respect to a fixed person in the group. In this case, the goal is to determine the probability that at least one of the other people has the same birthday as the chosen birthday person.

Unlike the classical birthday problem, in which all people in the group are compared to each other, the birthday person problem considers only comparisons involving a single individual. In this way, the number of comparisons performed is significantly lower, which produces a slower growth of coincidence probabilities.

\subsection{Formulation of the problem}

Given a fixed birthday, determine the probability that in a group with \(n\) people, at least one of the other people has the same birthday as the birthday person.

\subsection{Solution to the problem}

Initially, for the probabilistic modeling of the problem, we adopt the following hypotheses: a year with \(365\) days, in which all dates are equally likely; birthdays are independent; and February \(29\) is disregarded.

Let the event be

\[
A=\{\text{at least one of the other people has the same date as the birthday person}\}.
\]

We want to calculate the probability of occurrence of the event \(A\), that is,

\[
P(A).
\]

To simplify the calculations, we initially consider the complementary event of \(A\), defined by

\[
A^c=\{\text{none of the other people have the same birthday as the birthday person}\}.
\]

Observe that, since \(A\) and \(A^c\) are complementary events in \(\Omega\),

\[
A\cup A^c=\Omega.
\]

Necessarily, one of the two events occurs in our sample space: either there is a birthday coincidence, or all birthday dates are distinct. Furthermore,

\[
A\cap A^c=\emptyset,
\]

since an event and its complement cannot occur simultaneously. Since \(A\) and \(A^c\) are disjoint, by the additivity axiom,

\[
P(A\cup A^c)=P(A)+P(A^c).
\]

But

\[
A\cup A^c=\Omega.
\]

We obtain

\[
P(\Omega)=P(A)+P(A^c).
\]

By Kolmogorov's first axiom,

\[
P(\Omega)=1.
\]

Substituting \(P(\Omega)=1\), we obtain

\[
1=P(A)+P(A^c).
\]

Now, isolating \(P(A)\), one subtracts \(P(A^c)\) from both sides of the equality,

\[
1-P(A^c)=P(A)+P(A^c)-P(A^c).
\]

Since

\[
P(A^c)-P(A^c)=0,
\]

it follows that

\[
1-P(A^c)=P(A).
\]

Therefore,

\[
P(A)=1-P(A^c).
\]

Thus, it is enough to calculate \(P(A^c)\), that is, to determine the probability that none of the other people has the same birthday as the birthday person.

Now we calculate this probability. For this, we define the following events:

\begin{itemize}
    \item \(E_1\): the first person has a different birthday than the birthday person;
    \item \(E_2\): the second person has a different birthday than the birthday person;
    \item \(E_3\): the third person has a different birthday than the birthday person;
    \item \(\cdots\)
    \item \(E_{n-1}\): the \((n-1)\)-th person has a different birthday than the birthday person.
\end{itemize}

So, under these conditions,

\[
A^c=E_1\cap E_2\cap E_3\cap\cdots\cap E_{n-1}.
\]

Since birthdays are independent of each other, the events \(E_1,E_2,E_3,\ldots,E_{n-1}\) are also independent. By the definition of independence,

\[
P(A^c)
=
P(E_1)\,
P(E_2)\,
P(E_3)
\cdots
P(E_{n-1}).
\]

For each person, there are \(364\) favorable dates among the \(365\) days of the year, since only the birthday of the fixed person must be avoided. Thus,

\[
P(E_1)=\frac{364}{365},\qquad
P(E_2)=\frac{364}{365},\qquad
P(E_3)=\frac{364}{365}.
\]

Continuing the process, we obtain

\[
P(A^c)
=
\frac{364}{365}\cdot
\frac{364}{365}\cdot
\frac{364}{365}
\cdots
\frac{364}{365}.
\]

Since the same factor appears repeatedly \(n-1\) times, it follows that

\begin{equation}
P(A^c)
=
\left(\frac{364}{365}\right)^{n-1}.
\label{eq:aniversariante}
\end{equation}

We know that

\[
P(A)=1-P(A^c).
\]

Substituting Equation \ref{eq:aniversariante}, we obtain

\[
P(A)
=
1-\left(\frac{364}{365}\right)^{n-1}.
\]

Therefore, the probability that at least one person has the same birthday as the birthday person is

\begin{equation}
P(A)
=
1-\left(\frac{364}{365}\right)^{n-1}.
\label{eq:prob_aniversariante}
\end{equation}

\subsection{Example 2}

Consider a group with \(n=23\) people, including the previously fixed birthday person. We wish to determine the probability that at least one of the other people has the same birthday as that birthday person. Thus, by Equation \ref{eq:prob_aniversariante}, the probability is given by

\begin{equation}
P(A)
=
1-\left(\frac{364}{365}\right)^{22}.
\label{eq:exemplo_aniversariante}
\end{equation}

Applying the power property of a fraction, we obtain

\[
P(A)=1-\frac{364^{22}}{365^{22}}.
\]

It is noted that

\[
\frac{364}{365}\approx0{,}99726027.
\]

Substituting this value into Equation \eqref{eq:exemplo_aniversariante}, we obtain

\[
P(A)=1-(0{,}99726027)^{22}.
\]

Calculating the power,

\[
(0{,}99726027)^{22}\approx0{,}94143.
\]

Now, performing the subtraction,

\[
P(A)=1-0{,}94143.
\]

Therefore,

\[
P(A)\approx 0{,}05857.
\]

In percentage terms,

\[
P(A)\approx 5{,}857\%.
\]

The result indicates that, in a group with \(23\) people, including the previously fixed birthday person, the probability that at least one of the other people has the same birthday as that birthday person is approximately \(5{,}857\%\).

\subsection{Comparison with the birthday problem}
\label{subsec:comparacao}

It is observed that the growth of the probabilities in the birthday person problem occurs more slowly when compared to the classical birthday problem. This occurs because, while the classical model considers all possible comparisons between pairs of individuals, the birthday person problem restricts the comparisons to only one fixed person in the group.

In the birthday problem, each individual can be compared with everyone else, causing the number of comparisons to increase quadratically as a function of the size of the group. In contrast, in the birthday person problem, only comparisons involving the previously chosen birthday person are considered. Thus, for a group with \(n\) people, only \(n-1\) comparisons are performed.

\begin{figure}[H]
\centering
\begin{tikzpicture}
\begin{axis}[
    width=13cm, height=7cm,
    xlabel={$n$ (number of people)},
    ylabel={$P(A)$},
    xmin=0, xmax=80,
    ymin=0, ymax=1.05,
    xtick={0,10,...,80},
    ytick={0,0.2,...,1},
    grid=major,
    legend pos=north west,
    legend style={font=\small},
]
\addplot[blue, thick, mark=*] coordinates {
    (1, 0.000) (5, 0.027) (10, 0.117) (15, 0.253) (20, 0.411)
    (23, 0.507) (25, 0.569) (30, 0.706) (35, 0.814) (40, 0.891)
    (50, 0.970) (60, 0.994) (70, 0.999) (80, 1.000)
};
\addlegendentry{Birthday problem}
\addplot[red, thick, mark=square] coordinates {
    (1, 0.000) (5, 0.011) (10, 0.024) (15, 0.038) (20, 0.051)
    (23, 0.059) (25, 0.064) (30, 0.077) (35, 0.089) (40, 0.102)
    (50, 0.126) (60, 0.150) (70, 0.173) (80, 0.196)
};
\addlegendentry{Birthday person problem}
\end{axis}
\end{tikzpicture}
\caption{Comparison between the birthday problem and the birthday person problem: the latter grows much more slowly because it involves only \(n-1\) comparisons with a single fixed person.}
\label{fig:comparacao}
\end{figure}
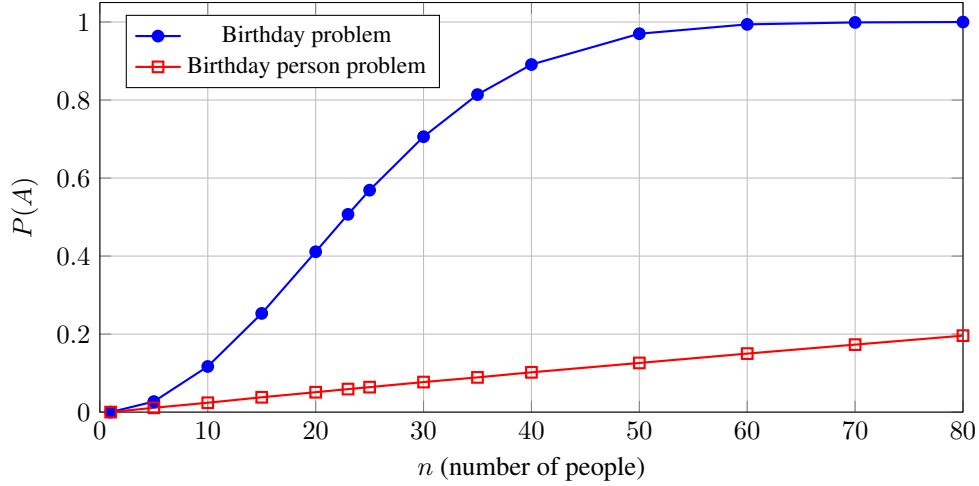

This difference explains why, in a group of \(23\) people, the probability that at least two of them share the same birthday is greater than \(50\%\), while the probability that one of them shares the same birthday as a previously fixed birthday person is only approximately \(5{,}857\%\). Figure \ref{fig:comparacao} and Table \ref{table:comparacao} summarize this contrast.

\begin{table}[H]
\centering
\begin{tabular}{|c|c|c|}
\hline
\textbf{Number of people} & \textbf{Birthday problem} & \textbf{Birthday person problem} \\
\hline
\(10\) & \(11{,}69\%\) & \(2{,}44\%\) \\
\(20\) & \(41{,}14\%\) & \(5{,}08\%\) \\
\(23\) & \(50{,}73\%\) & \(5{,}86\%\) \\
\(30\) & \(70{,}63\%\) & \(7{,}65\%\) \\
\(50\) & \(97{,}04\%\) & \(12{,}60\%\) \\
\(60\) & \(99{,}41\%\) & \(14{,}99\%\) \\
\hline
\end{tabular}
\caption{Numerical comparison between the birthday problem and the birthday person problem.}
\label{table:comparacao}
\end{table}

Therefore, although the two problems have a similar probabilistic foundation and use the principle of the complementary event to calculate the probabilities, their behaviors differ significantly due to the number of comparisons considered in each model. While the birthday problem involves all possible comparisons between pairs of individuals, the birthday person problem is limited to comparisons with a single person, producing a considerably slower growth in the probabilities of coincidence.

\section{Applications}
\label{sec:applications}

Although the birthday problem is frequently presented as a classical example of probability theory, its importance goes far beyond the mathematical context in which it was originally formulated. Currently, this problem has relevant applications in areas such as computing, information security, cryptography, and the analysis of large data sets, in which the occurrence of coincidences plays a fundamental role.

\subsection{Collisions in hash functions}

One of the best-known applications is related to hash functions, algorithms used to transform information of arbitrary size into numerical sequences of fixed size. Ideally, different information should produce different results. However, as the number of possible outputs is limited, there is the possibility that two different inputs generate the same result, a phenomenon known as a \emph{collision} \cite{stallings2015}.

The relationship between this phenomenon and the birthday paradox becomes evident when it is observed that, even if there are a large number of possible outcomes, the probability of coincidences occurring grows rapidly as the number of elements analyzed increases. In the context of birthday dates, the paradox shows that in a relatively small group of people, the chance that two of them share the same date of birth is already surprisingly high. Similarly, in hash functions, increasing the amount of information processed significantly increases the likelihood that two distinct inputs will produce the same result. Thus, the birthday problem provides an important probabilistic interpretation for understanding and estimating the occurrence of collisions in computational systems.

\begin{figure}[H]
    \centering
    \includegraphics[width=15cm]{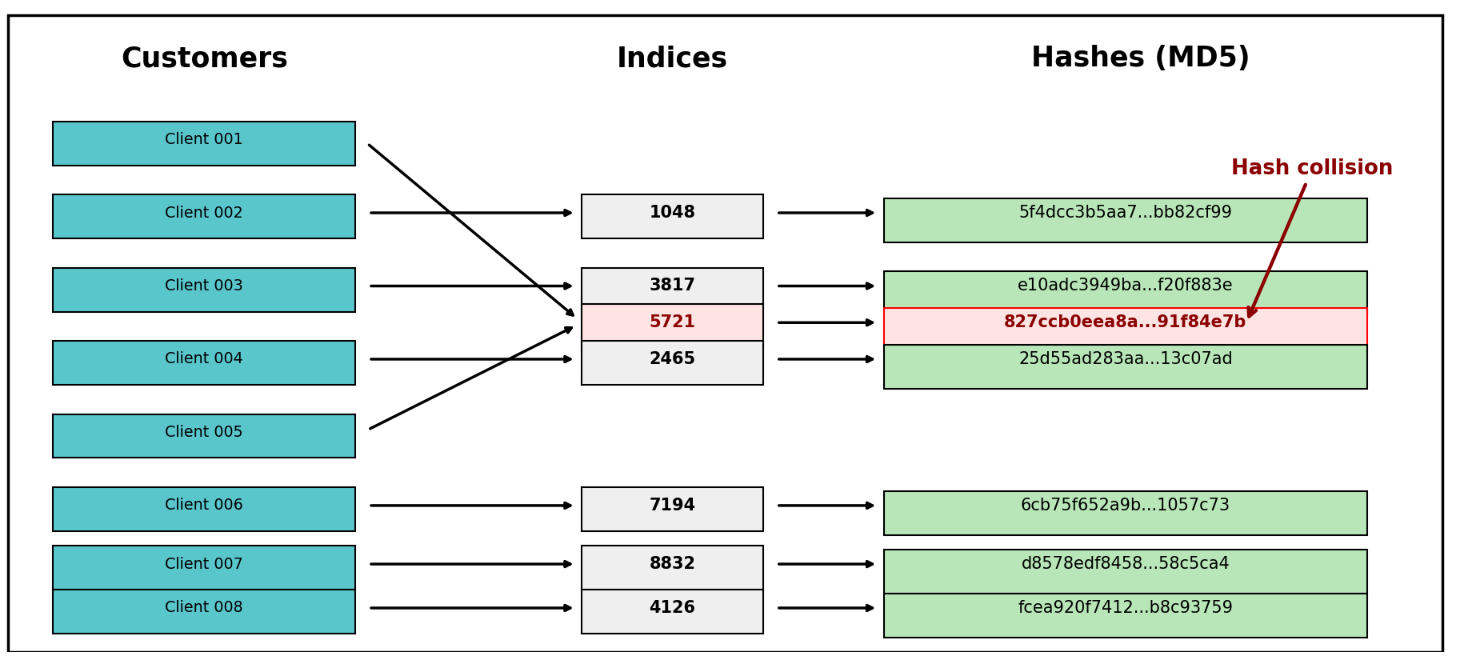}
    \caption{Schematic representation of a hash collision in a storage structure.}
    \label{fig:hash_armazenamento}
\end{figure}

\begin{figure}[H]
    \centering
    \includegraphics[width=15cm]{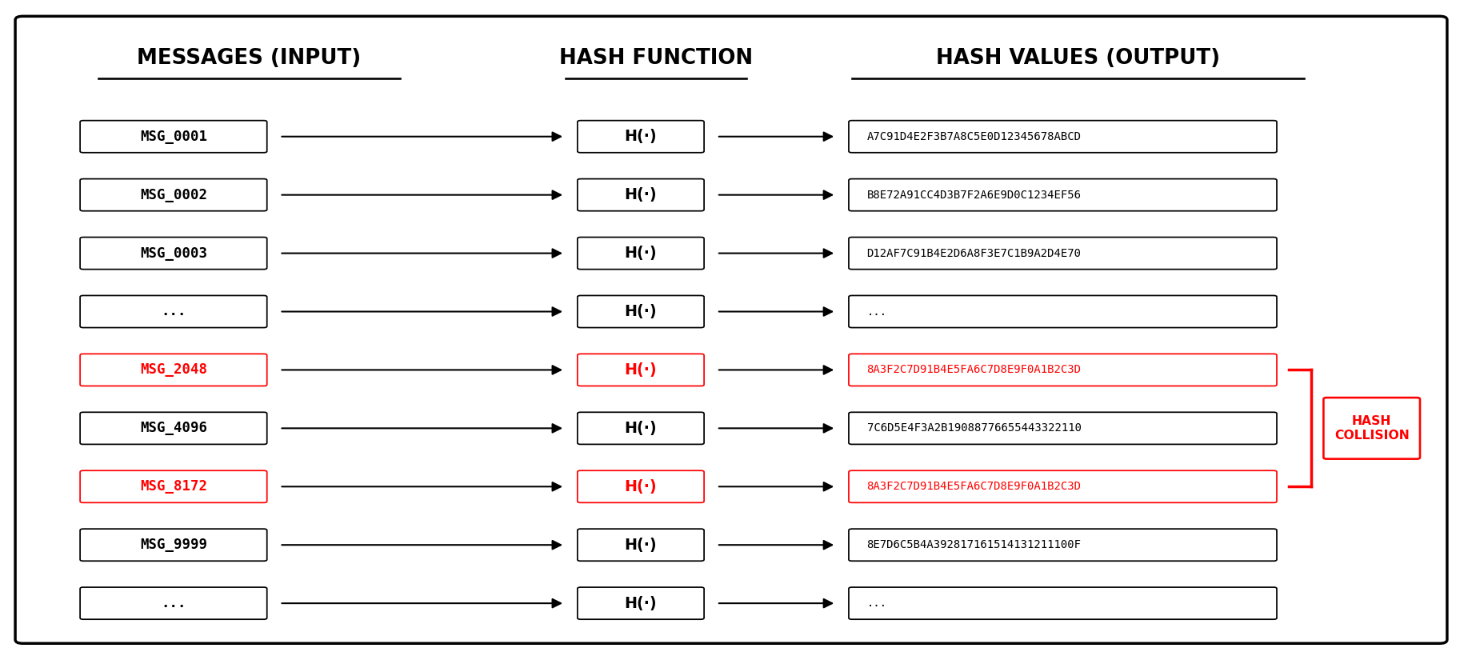}
    \caption{Schematic representation of a hash collision in cryptographic algorithms.}
    \label{fig:hash_criptografia}
\end{figure}

Figures \ref{fig:hash_armazenamento} and \ref{fig:hash_criptografia} present schematic representations of this process. In them, different information is associated with storage positions through a hash function, evidencing situations in which different entries end up being directed to the same index. These collisions illustrate, in a visual way, a probabilistic behavior similar to that observed in the birthday problem.

These representations were elaborated considering a set of independent inputs, represented by clients and messages, whose data is processed by a hash function and associated with storage indexes or hash values. To highlight the phenomenon discussed in this section, cases were highlighted in which different inputs produce the same output value, characterizing the occurrence of a collision.

In this representation, different inputs are processed by a hash function and associated with storage positions. The coincidence between two results evidences the occurrence of a collision, demonstrating that distinct inputs can be directed to the same index, even when the hash function is designed to distribute the data evenly. In this way, the illustration allows us to understand, intuitively, how this phenomenon can occur in computational applications and why its analysis is directly related to the birthday problem.

\subsection{The birthday attack}
\label{subsec:birthday_attack}

This interpretation has great relevance in the area of computer security, especially in the analysis of the reliability of cryptographic algorithms and digital authentication systems. In particular, the so-called \emph{birthday attack} relies on probabilistic properties similar to those of the birthday problem to exploit collisions in hash functions \cite{stallings2015}.

Mathematically, consider a hash function whose output has \(b\) bits, so that there are \(N=2^{b}\) possible hash values. Assuming that the hash behaves like a random function with uniformly distributed outputs, the probability that no collision occurs after evaluating \(n\) distinct inputs is, by analogy with Equation \eqref{eq:complementar},

\[
P(\text{no collision})
=
\prod_{k=0}^{n-1}\left(1-\frac{k}{N}\right)
\approx
\exp\left(-\frac{n(n-1)}{2N}\right).
\]

Consequently, the probability of at least one collision is

\[
P(\text{collision})
\approx
1-\exp\left(-\frac{n(n-1)}{2N}\right).
\]

Setting this probability equal to $1/2$ and solving for \(n\), one obtains the classical estimate

\begin{equation}
n_{1/2}\approx 1{,}1774\sqrt{N}=1{,}1774\cdot 2^{b/2}.
\label{eq:birthday_attack}
\end{equation}

This result has profound practical implications: the security of a hash function against collision attacks is limited not by \(2^{b}\), but by approximately \(2^{b/2}\). In other words, an attacker only needs to evaluate on the order of the square root of the output space to find a collision with non-negligible probability. Equation \eqref{eq:birthday_attack} is the mathematical core of the birthday attack.

\subsection{Numerical implications}

Table \ref{table:hash_attack} illustrates the practical consequences of Equation \eqref{eq:birthday_attack} for the hash functions most commonly used in cryptographic systems.

\begin{table}[H]
\centering
\begin{tabular}{|l|c|c|c|}
\hline
\textbf{Hash} & \textbf{Output ($b$ bits)} & \textbf{$2^{b}$ (space)} & \textbf{$n_{1/2}\approx 2^{b/2}$} \\
\hline
MD5     & $128$ & $2^{128}$ & $\approx 2^{64}\approx 1{,}8\times 10^{19}$ \\
SHA-1   & $160$ & $2^{160}$ & $\approx 2^{80}\approx 1{,}2\times 10^{24}$ \\
SHA-256 & $256$ & $2^{256}$ & $\approx 2^{128}\approx 3{,}4\times 10^{38}$ \\
\hline
\end{tabular}
\caption{Order of magnitude of the number of hash evaluations required to find a collision with probability $1/2$, according to Equation \eqref{eq:birthday_attack}.}
\label{table:hash_attack}
\end{table}

For MD5, although the output space has $2^{128}\approx 3{,}4\times 10^{38}$ values, a collision is expected after roughly $2^{64}\approx 1{,}8\times 10^{19}$ evaluations, which is within reach of modern computing clusters and has indeed been demonstrated in practice. For SHA-1, the corresponding bound is $2^{80}$, a value that, while still large, motivated the transition to SHA-2 and SHA-3 families. For SHA-256, the collision bound $2^{128}$ remains far beyond current computational capabilities, which is one of the reasons why it is widely recommended in security-critical applications. This discussion shows how the birthday problem, formulated in the 1930s as a purely combinatorial curiosity, has become a fundamental tool for assessing the security of modern cryptographic systems.

\subsection{Further applications and countermeasures}

In addition to cryptographic hash functions, the birthday problem has applications in areas such as load balancing in distributed hash tables (DHTs), duplicate detection in large data sets, and the analysis of pseudorandom number generators. In the context of information security, an important countermeasure against birthday attacks is the use of \emph{salted} hashes, in which a random value (salt) is concatenated to the input before hashing. Although salting does not increase the output size, it prevents precomputed collision tables and significantly raises the cost of practical attacks. Combined with larger digests (e.g., SHA-256 or SHA-3), salting provides robust protection against the class of attacks derived from the birthday problem.

From a mathematical point of view, the birthday problem remains an important tool for the study of applied probability and combinatorial analysis and is often used to introduce concepts such as complementary probabilities, combinatorial growth, and intuitive interpretation of random phenomena. In this way, it stands out not only as a classic result of probability theory, but also as a model capable of explaining and analyzing phenomena present in different areas of mathematics and computing.

\section{Conclusion}

The study developed made it possible to understand mathematically the behavior of the probabilities of coincidences between birthdays, highlighting results that frequently contradict intuitive perception. The analysis carried out showed that the occurrence of coincidences becomes much more probable than is normally imagined, even in relatively small groups of people.

The results obtained showed that the probability of coincidences increases rapidly as the number of people considered increases. In particular, it was observed that, in the birthday problem, relatively small groups already produce high probabilities of coincidence, while in the birthday person problem this growth occurs more slowly, due to the smaller number of comparisons involved.

One of the central aspects of the study was the use of complementary probability. In both problems, the calculation of the desired probability was performed in a simpler way through the initial determination of the probability of the complementary event, evidencing the importance of this technique in solving probabilistic problems.

Furthermore, the combinatorial interpretation proved to be fundamental for understanding the results obtained. The growth in the number of possible pairs among individuals explains why apparently unlikely coincidences become increasingly probable as the group size increases. The Poisson approximation \eqref{eq:poisson} provided a compact analytical expression for the probability curve, and the comparison with the birthday person problem clarified how the growth order (quadratic versus linear) determines the qualitative behavior of each model.

Finally, it was found that the concepts studied have applications that go beyond the theoretical context of mathematics, being present in areas such as computing, cryptography and information security, especially in the analysis of collisions in hash functions and in the so-called birthday attack. The mathematical formulation of the birthday attack (Equation \eqref{eq:birthday_attack}) and the numerical estimates presented in Table \ref{table:hash_attack} show that the birthday problem, far from being a mere curiosity, is a practical tool for evaluating the security of cryptographic systems. Thus, the birthday problem stands out as an important example of probability theory, capable of relating fundamental mathematical concepts with applications present in different scientific and technological contexts.


\end{document}